\documentclass[11pt,a4paper]{article}
\usepackage{a4wide}
\usepackage{theorem}
\usepackage{amsmath,amssymb,calc}
\usepackage[english]{babel}
\usepackage{graphicx,array}
\usepackage[numbers]{natbib}
\usepackage{tikz}
\usepackage{color} %Mag later weg, geeft voor mij aan waar ik nog even naar moet kijken
\usepackage{enumerate}
\usepackage[toc,page]{appendix}
\usepackage[charter]{mathdesign}
\usepackage{algorithm}
\usepackage{algpseudocode}
\usepackage{tikz}

\graphicspath{{./},{./images/}}

\newtheorem{theorem}{Theorem}%[section]
\newtheorem{assumption}[theorem]{Assumption}

\newtheorem{remark}[theorem]{Remark}
\newtheorem{proposition}[theorem]{Proposition}

\numberwithin{equation}{section}
\providecommand{\href}[2]{#2}

\newcommand{\beq}{\begin{equation}}
\newcommand{\eq}{\end{equation}}

\newcommand{\eqan}[1]{\begin{align} #1 \end{align}}
\newcommand{\e}{{\rm e}}
\renewcommand{\d}{{\rm d}}

\newcommand\blfootnote[1]{%
  \begingroup
  \renewcommand\thefootnote{}\footnote{#1}%
  \addtocounter{footnote}{-1}%
  \endgroup
}

\newcommand{\uu}{u}
\newcommand{\zz}{z}
\newcommand{\rr}{z}
\newcommand{\pgfA}{A}

\newcommand{\FF}{\mathcal{F}}
\newcommand{\ii}{{\rm i}}

{}
\newcommand{\ch}[1]{}
\newcommand{\chr}[1]{}

\title{Spitzer's identity for discrete random walks\blfootnote{Department of Mathematics and Computer Science, Eindhoven University of Technology, P.O. Box 513, 5600MB Eindhoven, The Netherlands. Email: \href{mailto:a.j.e.m.janssen@tue.nl}{a.j.e.m.janssen@tue.nl} and \href{mailto:j.s.h.v.leeuwaarden@tue.nl}{j.s.h.v.leeuwaarden@tue.nl}}
}

\author{A.J.E.M. Janssen \and J.S.H. van Leeuwaarden}

\begin{document}
\maketitle

\begin{abstract}
Spitzer's identity describes the position of a reflected random walk over time in terms of a bivariate transform. Among its many applications in probability theory are congestion levels in queues and random walkers in physics. We present a new derivation of Spitzer's identity under the assumption that the increments of the random walk have bounded jumps to the left. This mild assumption facilitates a proof of Spitzer's identity that only uses basic properties of analytic functions and contour integration. The main novelty, believed to be of broader interest, is a reversed approach that recognizes a factored polynomial expression as the outcome of Cauchy's formula. 
%Most of the research done on the FCTL queue assumes that the vehicles arrive at the intersection according to a Poisson process and focusses on deriving formulas for the mean queue length at the end of green periods and the mean delay. For a class of discrete arrival processes, including the Poisson process, we derive the probability generating function of both the queue length and delay, from which the whole queue length and delay distribution can be obtained.
%This allows for the evaluation of performance characteristics other than the mean, such as the variance and percentiles of the distribution.
%
%We discuss the numerical procedures that are required to obtain the performance characteristics, and give several numerical examples.

%\bigskip\noindent\textbf{Keywords:} fluctuation theory; transform methods; complex analysis\\
%\noindent {\bfseries AMS 2010 Subject Classification}. 60E10, 60J10, 60K25, 68M20, 90B20

\end{abstract}

%\tableofcontents

\section{Introduction}

Random walks are ubiquitous in the modern stochastics literature. This paper deals with the one-dimensional random walk on $\mathbb{Z}$ describing the partial sums $S_n:=X_1+\cdots+ X_n$ ($S_0:=0$) of i.i.d.~random variables $X_1,X_2,\ldots$. The stochastic process $(S_n,n\geq 0)$, the random walk with steps $X_n$, arises in many areas of science to describe the evolution of certain objects subject to random fluctuations, including random walkers in physics, congestion levels in queueing theory and capital positions in insurance mathematics \cite{comtet2005precise,janssen2007lerch}. If we indeed think of a random walk $(S_n,n\geq 0)$ as modeling capital or congestion, the monetary/physical interpretation means that large values of $S_n$ are of particular interest, and it is natural to consider the sequence $\bar M_n:=\max\{S_0,S_1,\ldots,S_n\}$. The study of $(\bar M_n,n\geq 0)$ and related quantities is referred to as the
fluctuation theory of the random walk, a central topic in most classic probability textbooks \cite{asmussen2003,chung2001course,cohen1982,feller1971,prabhu1998}.

Fluctuation theory became highly topical by the rise of queueing theory in the first half of the twentieth century, with foundational works of A.K.~Erlang and F.~Pollaczek (see the historical account in \cite{janssen2008back}) and as primary example the waiting time process in the single-server queue. Let consecutive customers arriving to a single-server queue be numbered $n=1,2,\ldots$. Denote by $B_{n}$ the service time of customer $n$, and by $C_{n}$ the time between the arrivals of customers $n$ and $n+1$. Then with $X_{n+1}:=B_n-C_n$, and $W_n$ the waiting time of customer $n$, the waiting time process (also known as the Lindley process) is given by
\eqan{\label{lindleyprocess}
W_{n+1}=(W_n+X_{n+1})^+,  \quad n=0,1,\ldots,
}
with $x^+=\max(0,x)$, $W_0$ assumed zero, and $X_1,X_2,\ldots$ i.i.d.~random variables.
It is readily seen that $W_{n}$  equals
 $M_n=\max\{0,X_n,X_n+X_{n-1},\ldots,X_n+\ldots+X_1\}$, and hence $W_n$, $\bar M_{n}$ and $M_n$ are all equal in distribution.
This connects queueing theory to many seemingly unrelated questions in applied probability, combinatorics and physics about $M_n$ (see e.g.~\cite{asmussen2003,banderierflajolet2002,feller1971,lalley2001}).

%The following crucial observation was made by Lindley \cite{lindley1952}. Consider
%a one-dimensional random walk of a particle starting from the origin and moving in
%jumps at discrete time intervals. Let $W_n$ be the distance from the origin after
%the $n$th jump of size $X_n$. Suppose further that there is an impenetrable barrier at the
%origin, so that the walk is never allowed to enter the negative half of the real line but
%is held at the origin whenever a jump would have carried it over the barrier, and does
%not move again until the value of $W_n$ is again positive. Then this one-dimensional {\it reflected} random walk is clearly described by the Lindley process in \eqref{lindleyprocess}.

For general distributions of $X_n$, pioneering work of Pollaczek \cite{pollaczek1930b,pollaczek1930a} (see also \cite[Section~5]{kingman1962}) resulted in formal solutions of the distribution of $M_n$ in terms of complex contour integrals (see \cite{Abate1993,boon2017pollaczek,Janssen2015} for the algorithmic aspects of these contour integrals). Another approach was taken by Spitzer \cite{spitzer1956}, who used combinatorial arguments to establish the identity now bearing his name.
%\begin{align}
%\FF(\uu,\zz )=\sum_{n\geq 0}\uu^n\mathbb{E}(\zz ^{W_n}).
%\end{align}
\begin{theorem}\label{thmgeneral} {\rm(Spitzer's identity \cite{spitzer1956})} For $|u|<1$ and $t\in \mathbb{R}$,
\begin{align}\label{ssssp}
\sum_{n\geq 0}\uu^n\mathbb{E}(\e ^{\ii t M_n})
&=\exp\Big\{\sum_{l=1}^{\infty}\frac{u^l}{l}\mathbb{E}(\e^{\ii t S_{l}^+})\Big\}.
\end{align}
\end{theorem}

In the course of a century, Spitzer's identity has been proved by several methods, some combinatorial, and some analytical.
The most common method in the literature to derive Spitzer's identity is to transform the Lindley process \eqref{lindleyprocess} into an integral equation for the distribution function of consecutive queue lengths, which is of the Wiener-Hopf type. Switching then from distribution functions to characteristic functions gives rise to a functional equation which is amenable to Wiener-Hopf factorization. Wiener-Hopf factorization was first applied to the stationary single-server queue by Smith \cite{smith1953}. Being a corner stone of applied probability, Spitzer's identity is derived in the famous textbooks \cite[Chapter 8]{asmussen2003}, \cite[Chapter 8]{chung2001course}, \cite[Chapter 5]{cohen1982}, \cite[Chapter 17]{feller1971} and \cite[Chapter 1]{prabhu1998}, and in all these books the method of choice is the Wiener-Hopf factorization.

In an expository paper on the single-server queue, Kingman \cite{kingman1962} reviews the analytical and combinatorial methods to establish Spitzer's identity. According to Kingman \cite{kingman1962}, the proof methods for the general single-server queue leading to Spitzer's identity used by Pollaczek (contour integrals), Smith (Wiener-Hopf), and Spitzer (combinatorics) turn out to be all variants of a heavily disguised algebraic structure, essentially having to do with the special projection properties of the $\max(0,\cdot)$ operator. This observation of Kingman, however, applies to the situation when a solution is required for quite general distributions of interarrival and service times, while much of the attention of queueing theorists has been directed to problems in which one or both of these distributions obey additional assumptions, so that more specialized techniques are applicable. These more restricted cases include random walks that are skip-free, so taking steps of maximal size one, in the positive (GI/M/1) or negative (M/G/1) direction, and the case in which the random variables $X_n$ are integer-valued (or have a lattice distribution, that is, concentrated on the integer multiples of some real number). Using skip-free or lattice properties allows to analyze the process \eqref{lindleyprocess} by methods that do not rely on Spitzer's identity \eqref{ssssp}.

Indeed, using lattice assumptions, one can consider generating functions (or Laplace transforms) instead of characteristic functions, and the analytic properties of generating functions make it possible to solve for the
bivariate generating function
 \begin{align}\label{}
\FF(\uu,\zz )=\sum_{n\geq 0}\uu^n\mathbb{E}(\zz ^{M_n})
\end{align}
using the so-called kernel method. This method was first applied by Crommelin \cite{crommelin1932,crommelin1934} to the stationary version of Lindley's equation \eqref{lindleyprocess} with $X_n=A_n-s$, $A_n$ Poisson distributed and $s$ a nonnegative integer. The idea is to find the generating function by solving the functional equation for $\FF(u,z)$ using a factorization in terms of the complex-valued roots of the kernel $z^s-u\mathbb{E}(z^{A_n})$ (see \eqref{10} below). The function $\FF(u,z)$ is also called  Green's function, and like \eqref{ssssp} completely characterizes the  distribution of the position of the reflected random walk at all points in time. For random walks with lattice increments $X_n$, the kernel method has been applied to obtain $\FF(u,z)$  in \cite{banderierflajolet2002,lalley2001}. Building on this general result for $\FF(u,z)$, we are able to construct a fully analytic proof of the following version of Spitzer's identity.
\begin{proposition}\label{thmdiscrete} {\rm(Spitzer's identity for discrete queues)}
Consider the stochastic process $(M_n,n\geq 0)$ with $X_n$ an integer-valued random variable with a support that is bounded from below.
Then, for $|u|<1$ and $|z|\leq 1$,
\begin{align}\label{sp1}
\FF(\uu,\zz )
&=\exp\Big\{\sum_{l=1}^{\infty}\frac{u^l}{l}\mathbb{E}(\zz^{S_{l}^+})\Big\}.
\end{align}
\end{proposition}
Let us stress that Proposition \ref{thmdiscrete} is less general than Theorem \ref{thmgeneral}.
Over the years the identity has appeared under various conditions. Pollaczek \cite{pollaczek1930b,pollaczek1930a} assumed light-tailed distributions,
%needed the condition that \chr{$\mathbb{E}(\e^{-\theta X_n})<\infty$ for some $\theta<0$}, related to the assumption of light-tailed distributions.
and Spitzer's combinatorial proof \cite{spitzer1956} required the assumption of integer-valued increments, like in our setting. Building on the initial results of Pollaczek, Spitzer, and others, it was realized that a certain care was needed to give rigorous proofs when stretching the assumptions on the increment distribution, with Theorem \ref{thmgeneral} as the general version in which it appears in \cite{asmussen2003,feller1971,prabhu1998}. For instance, one could ask whether the Laplace transform is obtained by replacing ``$\ii t$'' in \eqref{ssssp} by $-\theta$ with $\theta>0$. Or in our case, whether the generating function in Proposition \ref{thmdiscrete} can lead to a result for a Laplace transform, or even a characteristic function.

In Kingman's \cite{kingman1962} words, Pollaczek's proof constitutes a dazzling series of manipulations, using an elaborate machinery of contour integration. This also explains why many contemporaries considered Pollaczek's work rather impenetrable.  By using Laplace transforms instead of characteristic functions, Pollaczek was able to give an analytic derivation of Spitzer's identity. Like Pollaczek, our proof of Spitzer's identity relies strongly on analytic functions and contour integration. However, the combination of the additional assumptions of integer-valued and {\it bounded} increments, and a new idea that exploits the kernel method and the analytic properties of generating functions, leads to a short and transparent proof. Providing this proof for Proposition \ref{thmdiscrete}  is the principal goal of this paper. It would be relatively easy to extend the scope of Proposition \ref{thmdiscrete}.  Any continuous distribution can be approximated arbitrarily closely by a lattice distribution (see \cite{konheim1975,stadje1997}), and then the discrete version of Spitzer's identity in Proposition \ref{thmdiscrete} can be shown to hold under more general conditions on the distribution of the increments $X_n$ using continuity arguments. This would introduce quite a few additional technicalities, though.

\section{Proof of Spitzer's identity}
We shall now prove Proposition \ref{thmdiscrete}. We consider the discrete queue described by the Markov chain $(W_n,n\geq 0)$ in \eqref{lindleyprocess}
with $X_1,X_2,\ldots$ a sequence of i.i.d.~discrete random variables. We assume that  $X_n$ is integer valued, so that the random walk described by the process $(M_n,n\geq 0)$ lives on the set of nonnegative integers.
Let $X$ denote a generic random variable equal in distribution to $X_n$, and assume
\begin{equation}
X\in\{-s,-s+1,\ldots,-1,0,1,\ldots \}
\end{equation}
with $\mathbb{P}(X=-s)\in(0,1)$ and $s$ some positive integer. This means that we can write the probability generating function (pgf)
of $X$ as
$$
\mathbb{E}(\zz^X)=\pgfA(\zz)\zz^{-s},
$$
where $\pgfA(\zz)$ is the pgf of a nonnegative integer-valued random variable $A$, i.e.,
\begin{equation}\label{lindely2}
M_{n+1}=(M_n+A_{n+1}-s)^+, \quad n=0,1,\ldots
\end{equation}
with $A_1,A_2,\ldots$ a sequence of i.i.d.~discrete random variables with $A\stackrel{d}{=} A_1$. We assume that $A(z)$ is analytic in a disk $|z|<R$ with $R>1$.

The outline of the proof of
Proposition \ref{thmdiscrete} is now as follows.
We shall first obtain the following product representation for $\FF(\uu,\zz )$, which only requires elementary manipulations of the recursion relation \eqref{lindely2}.
\begin{proposition}\label{lemmafac} {\rm(Product representation)} For $|u|<1$ and $|z|\leq 1$,
\begin{equation}\label{10}
\FF(\uu,\zz )=\frac{1}{\zz^s-\uu \pgfA(\zz)}\prod_{k=0}^{s-1}\frac{\zz -\rr_k(\uu)}{1-\rr_k(\uu)}
\end{equation}
with $\rr_k(\uu)$, $k=0,\ldots, s-1$ the $s$ roots of $\zz^s-\uu\pgfA(\zz)$ within the unit disk $|z|\leq 1$.
 %denoted by $\rr_0(\uu),\rr_1(\uu),\ldots,\rr_{s-1}(\uu)$.
\end{proposition}
Proposition \ref{lemmafac} is well known, see e.g.~\cite{banderierflajolet2002,bruneelkim,lalley2001}, but we shall provide a concise derivation in Subsection~\ref{sub1} to make this paper self-contained. Continuing then with the representation
\eqref{10}, the next step in our proof of Spitzer's identity is to transform the expression
\eqref{10} into a contour-integral representation. This is the main idea in this paper: interpret \eqref{10} as the outcome of Cauchy's residue theorem, the classical tool from complex analysis to evaluate integrals of analytic functions along closed
curves, leading to the following result:
\begin{proposition}\label{thmpol} {\rm(Pollaczek integral)}
Let $0<v<1$. There is a $b\in(1,R)$ such that
\begin{equation}\label{polk}
\FF(u,z)=\frac{1}{1-u}\exp\Big(\frac{1}{2\pi i} \oint_{|w|=b} \frac{z-1}{(w-1)(w-z)}\ln\left(1-u w^{-s}A(w)\right)\d w\Big)
\end{equation}
holds for $|u|\leq v$ and $|z|\leq b$.
\end{proposition}
The derivation of Proposition \ref{thmpol} is given in Subsection~\ref{sub2}.
The proof of Spitzer's identity is then completed in Subsection~\ref{sub3} by series expansion of the functions in \eqref{polk} and
identifying the resulting contour integrals as probabilities using the inversion formula for generating functions.

\subsection{Product representation (Proof of Proposition~\ref{lemmafac})}\label{sub1}
Let $|z|\leq 1$. Observe that
\begin{align}
\mathbb{E}(\zz ^{M_{n+1}})&=\mathbb{E}(\zz ^{(M_n+X_{n+1})^+})\nonumber\\
&=\mathbb{E}(\zz ^{(M_n+X_{n+1})^+}{\bf 1}_{\left\{M_n+X_{n+1}\geq 0\right\}})+\mathbb{E}(\zz ^{(M_n+X_{n+1})^+}{\bf 1}_{\left\{M_n+X_{n+1}< 0\right\}})\nonumber\\
&=\mathbb{E}(\zz ^{M_n+X_{n+1}})-\mathbb{E}(\zz ^{M_n+X_{n+1}}{\bf 1}_{\left\{M_n+X_{n+1}< 0\right\}})+\mathbb{P}(M_n+X_{n+1}< 0).
\end{align}
Now $\mathbb{E}(\zz ^{M_n+X_{n+1}})=\mathbb{E}(\zz ^{M_n})\pgfA(\zz)\zz^{-s}$,
\begin{align}
\mathbb{E}(\zz ^{M_n+X_{n+1}}{\bf 1}_{\left\{M_n+X_{n+1}< 0\right\}})
&=\sum_{r=0}^{s-1}\mathbb{P}(M_n+A_{n+1}=r) \zz ^{r-s}
\end{align}
and
$
\mathbb{P}(M_n+X_{n+1}< 0)=\sum_{r=0}^{s-1}\mathbb{P}(M_n+A_{n+1}=r).
$
This gives
\begin{align}\label{22222}
\mathbb{E}(\zz ^{M_{n+1}})&=\mathbb{E}(\zz ^{M_n})\pgfA(\zz)\zz^{-s}+\sum_{r=0}^{s-1}\mathbb{P}(M_n+A_{n+1}=r)(1-\zz ^{r-s}).
\end{align}
We now derive an expression for the bivariate generating function
\begin{equation}
\FF(\uu,\zz )=\sum_{n\geq 0}\uu^n\mathbb{E}(\zz ^{M_n}),\quad |u|<1, \ |z|\leq 1,
\end{equation}
using a similar approach as in e.g.~\cite{banderierflajolet2002,bruneelkim,lalley2001}.
From \eqref{22222} we get
\begin{align}
\FF(\uu,\zz )&=\mathbb{E}(\zz ^{M_0})+\uu \FF(\uu,\zz )\pgfA(\zz)\zz^{-s}+\uu\sum_{r=0}^{s-1}(1-\zz ^{r-s})F_r(\uu)
\end{align}
with $F_r(\uu)=\sum_{n=0}^\infty \mathbb{P}(M_n+A_{n+1}=r) \uu^n$.
Upon some rewriting we arrive at
\begin{equation}\label{mainpgfqueue}
\FF(\uu,\zz )=\frac{N(\uu,\zz )}{\zz^s-\uu\pgfA(\zz)},
\end{equation}
where
\begin{equation}\label{2110}
N(\uu,\zz )=\zz^s\mathbb{E}(\zz ^{M_0})+\uu\sum_{r=0}^{s-1}(\zz^s-\zz ^{r})F_r(\uu).
\end{equation}
The numerator $N(\uu,\zz )$ is a polynomial in $\zz$ of degree $s$, with coefficients depending on $\uu$.
From $|z^s|=1>|u|=|uA(1)|\geq |uA(z)|$ for $|z|=1$ and $|u|<1$, we have by Rouch\'{e}'s theorem \cite{adan2006application} that the denominator $\zz^s-\uu\pgfA(\zz)$ has exactly $s$ zeros within the open unit disk denoted by $\rr_0(\uu),\rr_1(\uu),\ldots,\rr_{s-1}(\uu)$. The function $\FF(\uu,\zz )$ is analytic in the polydisk $|\uu|<1$, $|\zz|<1$. Therefore, the zeros in $|z|<1$ of the denominator in \eqref{mainpgfqueue} should also be the zeros of the numerator.

Hence,
\begin{equation}\label{eq7}
N(\uu,\zz )=\gamma(\uu)\prod_{k=0}^{s-1}(\zz -\rr_k(\uu)),
\end{equation}
where $\gamma(\uu)$ follows from \eqref{2110} and $N(\uu,1)=1$ since $W_0=0$.
%\begin{equation}\label{eq8}
%h(1,u)\ = \ \frac{1}{1-u} \ =  \ \frac{\gamma(u)\prod_{k=0}^{s-1}(1-\uu_k(u))}{1-u}.
%\end{equation}
We thus arrive at the expression \eqref{10}.
%Denote by $Q$ the limit of $W_n$ as $n$ goes to infinity (which exists if $\rho=\pgfA'(1)/s<1$), and let $\FF(\zz)=\mathbb{E}(\zz^Q)$.
% \begin{align}
%\FF(\zz)&=\lim_{\uu\rightarrow 1}(1-\uu)\FF(\uu,\zz ) \nonumber \\&= \lim_{\uu\rightarrow 1}\frac{1-\uu}{1-\rr_0(\uu)}\frac{\zz-\rr_0(\uu)}{\zz^s-\uu \pgfA(\zz)}\prod_{k=1}^{s-1}\frac{\zz -\rr_k(\uu)}{1-\rr_k(\uu)} \nonumber\\
%&=\frac{(s-\pgfA'(1))(\zz-1)}{\zz^s-\pgfA(\zz)}\prod_{k=1}^{s-1}\frac{\zz -\rr_k}{1-\rr_k},\label{genform}
%\end{align}
%with $\rr_k\equiv \rr_k(1)$ and where we have used that $\rr_0'(1)=1/(s-\pgfA'(1))$.

%\bibliographystyle{abbrv}

\subsection{Pollaczek integral (Proof of Proposition~\ref{thmpol})}\label{sub2}
Fix $v\in(0,1)$. Since %$A(w)$ is analytic in $|w|<R$ with $R>1$ and $A(1)=1$,
$|A(w)|\leq A(|w|)$ when $|w|<R$, and $A(1)=1$, there are points $a,b$ with $0<a<1<b<R$
such that $|w^s|>|uA(w)|$ holds for all $u$, $|u|\leq v$ and all $w$, $a\leq |w| \leq b$.
In particular, by Rouch\'{e}'s theorem, all $s$ zeros $z_k(u)$ of $w^s-uA(w)$ with $|w|\leq 1$ satisfy $|z_k(u)|<a$ while
$w^s-uA(w)$ is zero-free in the annulus $a\leq |w|\leq b$, provided that $|u|\leq v$.

Now fix $u\in(0,v)$ and $z\in(a,1)$. The function $\ln(\frac{z-w}{1-w})$ is analytic in $w\in\mathbb{C}\setminus[z,1]$ when we take the principal value $\ln$. Setting $k(w)=w^s-uA(w)$ the function $k'(w)/k(w)$ has its poles within $|w|=a$ at $w=z_k(u)$, $k=0,\ldots,s-1$. We allow here that several $z_k(u)$ coincide, in which case the residue of $k'(w)/k(w)$ at such a $w=z_k(u)$ equals the multiplicity of the zero of $k(w)$. By Cauchy's theorem we then have that
\begin{equation}\label{213}
\sum_{k=0}^{s-1}\ln\left(\frac{z-z_k(u)}{1-z_k(u)}\right)=\frac{1}{2\pi i} \oint_{|w|=a}\ln\left(\frac{z-w}{1-w}\right)\frac{k'(w)}{k(w)} \d w.
\end{equation}
Using that $k(w)$ is analytic and zero-free in $a\leq |w|\leq b$, we get, again by Cauchy's theorem,
\begin{align}\label{214}
\sum_{k=0}^{s-1}\ln\left(\frac{z-z_k(u)}{1-z_k(u)}\right)&=\frac{1}{2\pi i} \oint_{|w|=b}\ln\left(\frac{z-w}{1-w}\right)\frac{k'(w)}{k(w)} \d w \nonumber\\
&-\frac{1}{2\pi i} \oint_{\mathcal{C}}\ln\left(\frac{z-w}{1-w}\right)\frac{k'(w)}{k(w)} \d w,
\end{align}
where $\mathcal{C}$ is a contour contained in $a\leq |w|\leq b$ that encircles the branch cut $[z,1] $ once in positive sense.

\begin{figure}[t]
\begin{center}
\begin{tikzpicture}[scale=3]
\draw[blue, very thick](1,0) circle [radius=0.3];
\draw[blue, very thick](2,0) circle [radius=0.3];
\draw[white,fill=white](0.9,-0.02) rectangle (2.2,0.02);
\draw[blue, thick](1.293,-0.02) -- (1.707,-0.02);
\draw[blue, thick](1.293,0.02) -- (1.707,0.02);
\node[above right] at (1,0.25) {$C_z$};
\node[above right] at (2,0.25) {$C_1$};
%\draw[help lines](1,0) -- (1.21,0.21);
\draw[black](1,0) -- (1.21,0.21);
\node[right] at (1,0.17) {\small$\delta$};
%\draw[help lines](2,0) -- (2.21,0.21);
\draw[black](2,0) -- (2.21,0.21);
\node[right] at (2,0.17) {\small$\delta$};

\draw[black](0,0) -- (2.7,0);
%\draw[black](0,0) circle [radius=1];
%\draw[black](0,0) circle [radius=2.75];
\node[below] at (0,0) {$0$};
\node[below] at (1.0,0.00) {$z$};
\draw[black](0,0.03) -- (0,-0.03);
\draw[black](1,0.03) -- (1,-0.03);
\draw[black](2,0.03) -- (2,-0.03);
\draw[black](0.5,0.03) -- (0.5,-0.03);
\draw[black](2.5,0.03) -- (2.5,-0.03);
\node[below ] at (0.5,-0.0) {$a$};
\node[below ] at (2.5,-0.0) {$b$};
%\node[above] at (3.7,0) {$t_0$};
%\node[below] at (3.7,0) {$R_0$};
% C1 and Cw
\node[below] at (2,0) {$1$};
\draw[thick,blue,->] (1.51,0.02)--(1.5,0.02);
\draw[thick,blue,<-] (1.51,-0.02)--(1.5,-0.02);
\draw[thick,blue,->] (0.7,0.01)--(0.7,0);
\draw[thick,blue,->] (2.3,-0.01)--(2.3,0);
%\draw[black,fill=black](0.7,0.35) circle[radius=\pgflinewidth];
%\draw[black,fill=black](0.7,-0.35) circle[radius=\pgflinewidth];
%\draw[black,fill=black](0.4,0.5) circle[radius=0.01];
%\draw[black,fill=black](0.4,-0.5) circle[radius=0.01];
%\draw[ultra thick,black,dash pattern=on \pgflinewidth off 1.1cm] (0.4,0.5) to [out=160,in=90] (-0.6,0)  ;
%\draw[ultra thick,black,dash pattern=on \pgflinewidth off 1.1cm] (0.4,-0.5) to [out=200,in=-90] (-0.6,0);
%\draw[ultra thick,black,dash pattern=on \pgflinewidth off 1.03cm] (0.7,0.35) to [out=130,in=90] (-0.6,0)  ;
%\draw[ultra thick,black,dash pattern=on \pgflinewidth off 1.03cm] (0.7,-0.35) to [out=230,in=-90] (-0.6,0);
%\node[above] at (0.7,0.35) {$z_1$};
%\node[above] at (0.4,0.55) {$z_2$};
%\node[below right] at (0.6,-0.35) {$z_{g-1}$};
%\node[below right] at (0.3,-0.55) {$z_{g-2}$};
\node[above] at (1.5,0.02) {\small $L_+$};
\node[below] at (1.5,-0.02) {\small $L_-$};
\end{tikzpicture}
\end{center}
\caption{The four components $C_z$, $C_1$, $L_+$ and $L_-$ of the contour $\mathcal{C}$.}
\label{fig:contours}
\end{figure}
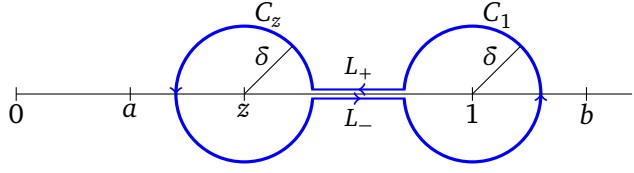

We shall now show that
\begin{align}\label{215}
\frac{1}{2\pi i} \oint_{\mathcal{C}}\ln\left(\frac{z-w}{1-w}\right)\frac{k'(w)}{k(w)} \d w =\ln(1-u)-\ln k(z).
\end{align}
We take for $\mathcal{C}$  the positively oriented contour consisting of the circles $C_z(\delta)$ and $C_1(\delta)$ of radii $\delta$ around $z$ and $1$, respectively, together with the line segments $L_{\pm}(\delta)=\{w = t \pm i 0 \,\,|\,\,z+\delta\leq t\leq 1-\delta\}$ where
 $0<\delta<\min\{z-a,\tfrac12 (1-z),b-1\}$. See Fig.~\ref{fig:contours} for an illustration of $\mathcal C$ and its four components.
For $t\in(z,1)$, we get (since $\ln$ is the principal value)
\begin{align}\label{216}
\ln\left(\frac{z-(t\pm i0)}{1-(t\pm i0)}\right)=\ln\left(\frac{t-z}{1-t}\right)\mp \pi i,
\end{align}
in which $\ln(\frac{t-z}{1-t})$ is real. Taking into account the orientation of $L_{\pm}$, this gives
\begin{align}\label{217}
&\frac{1}{2\pi i} \left(\int_{L_+}+\int_{L_-}\right)\ln\left(\frac{z-w}{1-w}\right)\frac{k'(w)}{k(w)} \d w \nonumber\\
&=\frac{1}{2\pi i} \int_{z+\delta}^{1-\delta}\left(-\left(\ln\left(\frac{t-z}{1-t}\right)-\pi i\right)+
\left(\ln\left(\frac{t-z}{1-t}\right)+\pi i\right)\right)\frac{k'(t)}{k(t)} \d t\nonumber\\
&=\int_{z+\delta}^{1-\delta}\frac{k'(t)}{k(t)}\d t=\ln(k(1-\delta))-\ln(k(z+\delta)),
%\ln(1-u)-\ln k(z).
\end{align}
where we observe that $k(t)=t^s-u A(t)>0$, $t\in[z,1]$. Using $k(1)=1-u$ and differentiability of $k(t)$, $z\leq t\leq 1$, we then get
\begin{align}\label{218}
&\frac{1}{2\pi i} \left(\int_{L_+}+\int_{L_-}\right)\ln\left(\frac{z-w}{1-w}\right)\frac{k'(w)}{k(w)} \d w =\ln(1-u)-\ln(k(z))+O(\delta).
%\ln(1-u)-\ln k(z).
\end{align}
Furthermore, since $\ln(\frac{z-w}{1-w})=O(\ln \delta)$, $w\in C_z(\delta)$ or $C_1(\delta)$,
\begin{align}\label{219}
&\frac{1}{2\pi i} \left(\oint_{C_z(\delta)}+\oint_{C_1(\delta)}\right)\ln\left(\frac{z-w}{1-w}\right)\frac{k'(w)}{k(w)} \d w =O(\delta\ln \delta).
%\ln(1-u)-\ln k(z).
\end{align}
Adding \eqref{218} and \eqref{219}, and letting $\delta\downarrow 0$, gives \eqref{215}.

Returning to \eqref{214}, we see that
\begin{align}\label{220}
\sum_{k=0}^{s-1}\ln\left(\frac{z-z_k(u)}{1-z_k(u)}\right)&=\frac{1}{2\pi i} \oint_{|w|=b}\ln\left(\frac{z-w}{1-w}\right)\frac{k'(w)}{k(w)} \d w
-\left(\ln(1-u)-\ln k(z)\right).
\end{align}
We rewrite the left-hand side of \eqref{220} as
\begin{align}\label{221}
\sum_{k=0}^{s-1}\ln\left(\frac{z-z_k(u)}{1-z_k(u)}\right)=\ln\left(\prod_{k=0}^{s-1}\frac{z-z_k(u)}{1-z_k(u)}\right)
\end{align}
with principal value logarithm at either side.
Here we use that for $a<z<1$ and $u>0$ the $z_k(u)$ come in conjugate pairs or are real in which case $z_k(u)<a<z<1$, causing the product at the right-hand side of \eqref{221} to be positive and the left-hand side of \eqref{221} to be real by cancellation of all imaginary parts of the principal logarithms. Then, from \eqref{10}, \eqref{220} and \eqref{221},
\begin{align}\label{222}
\ln\FF(\uu,\zz )&=-\ln(z^s-uA(z))+\ln\left(\prod_{k=0}^{s-1}\frac{z-z_k(u)}{1-z_k(u)}\right)\nonumber\\
&=-\ln(1-u)+\frac{1}{2\pi i} \oint_{|w|=b}\ln\left(\frac{z-w}{1-w}\right)\frac{\left(w^s-u A(w)\right)'}{w^s-u A(w)} \d w .
\end{align}
To establish \eqref{polk} from \eqref{222} for the $u$ and $z$ to which we have restricted ourselves, we compute
\begin{align}\label{223}
\frac{\left(w^s-u A(w)\right)'}{w^s-u A(w)}=\frac{s}{w}+\frac{\d }{\d w} \ln \left(1-u w^{-s}A(w)\right)
\end{align}
with analytic principal value logarithm since {$|u w^{-s}A(w)|<1$} when $0<u<v$, $a\leq |w|\leq b$. Now
\begin{equation}\label{224}
 \oint_{|w|=b} \frac{1}{w}\ln\left(\frac{z-w}{1-w}\right)\d w=0
\end{equation}
by Cauchy's theorem, where we replace $b$ in \eqref{224} by $c\geq b$ and let $c\to\infty$ while we observe that the integrand in \eqref{224} is $O(1/w^2)$.
Then
\begin{align}\label{225}
\ln\FF(\uu,\zz )
&=-\ln(1-u)+\frac{1}{2\pi i} \oint_{|w|=b}\ln\left(\frac{z-w}{1-w}\right)\frac{\d }{\d w} \ln \left(1-u w^{-s}A(w)\right) \d w .
\end{align}
Finally, by partial integration with continuous differentiable functions $\ln(\frac{z-w}{1-w})$  and
$\ln (1-u w^{-s}A(w)) $ on the closed contour $|w|=b$, we get
\begin{align}\label{226}
\ln\FF(\uu,\zz )
&=-\ln(1-u)+\frac{1}{2\pi i} \oint_{|w|=b}\left(\frac{\d }{\d w}\ln\left(\frac{z-w}{1-w}\right)\right) \ln \left(1-u w^{-s}A(w)\right) \d w ,
\end{align}
which establishes  \eqref{polk} for $0<u<v$, $a< z<1$. We get the validity of
\eqref{polk} on the full ranges $|u|\leq v$ and $|w|<b$ by analyticity, using that
$|u w^{-s}A(w)|<1$  when $|u|\leq v$ and $|w|<b$.

%\begin{equation}\label{212}
%\frac{k'(w)}{k(w)}=\frac{\left(w^s-uA(w)\right)'}{w^s-uA(w)}
%\end{equation}

\subsection{From integrals to infinite series}\label{sub3}
We continue from \eqref{226}, and develop for $|u|\leq v$ and $|z|<b=|w|$
\begin{equation}\label{227}
\frac{z-1}{(w-1)(w-z)}=\sum_{k=1}^\infty
\frac{1-z^k}{w^{k+1}}
\end{equation}
and
\begin{equation}\label{228}
 \ln \left(1-u w^{-s}A(w)\right)=-\sum_{l=1}^\infty\frac{1}{l}u^lw^{-sl}A^l(w).
\end{equation}
This gives by absolute and uniform convergence
\begin{align}\label{229}
\ln\FF(\uu,\zz )
&=-\ln(1-u)-\sum_{k=1}^\infty\sum_{l=1}^\infty \frac{1}{l}u^l(1-z^k)\frac{1}{2\pi i}\oint_{|w|=b}\frac{A^l(w)}{w^{k+sl+1}}\d w.
\end{align}
Observe that
\begin{align}\label{230}
\frac{1}{2\pi i}\oint_{|w|=b}\frac{A^l(w)}{w^{k+sl+1}}\d w=\mathbb{P}\Big(\sum_{i=1}^l A_i=k+sl\Big)=\mathbb{P}(S_l=k)
\end{align}
and use that $-\ln(1-u)=\sum_{l=1}^\infty u^l/l$ to get
\begin{align}\label{231}
\ln\FF(\uu,\zz )
&=\sum_{l=1}^\infty \frac{u^l}{l}-\sum_{l=1}^\infty \frac{u^l}{l}\sum_{k=1}^\infty(1-z^k) \mathbb{P}(S_l=k)\nonumber \\
&=\sum_{l=1}^\infty \frac{u^l}{l}\Big[(1-\mathbb{P}(S_l\geq 1))+\sum_{k=1}^\infty z^k \mathbb{P}(S_l=k)\Big],%\nonumber \\
%&=\sum_{l=1}^\infty \right]\nonumber \\
\end{align}
which gives \eqref{sp1}. The identity  \eqref{sp1} extends to all $|u|<1$ and $|z|\leq 1$ by analyticity and continuity, which completes the proof.

\section*{Acknowledgments}
%We thank the authors of \cite{ahmad} for sharing an early version of their work.
This work is supported by the NWO Gravitation Networks grant 024.002.003, an NWO TOP-GO grant and by an ERC Starting Grant.

\bibliographystyle{plain}
\bibliography{bibbook}

\end{document}